\documentclass[12pt]{elsarticle}
\usepackage{fullpage}
\usepackage{amsmath,amsfonts,amsthm}
\usepackage{upquote,graphicx,xcolor}
\usepackage{url}

\usepackage{pgfplots}
\pgfplotsset{compat=newest}

\usetikzlibrary{plotmarks}
\usetikzlibrary{arrows.meta}
\usepgfplotslibrary{patchplots}
\usepackage{grffile}

\newtheorem{definition}{Definition}[section]

\newcommand{\ee}{\text{e}}
\newcommand{\bu}{\mathbf u}
\newcommand{\bx}{\mathbf x}
\newcommand{\matlab}{\texttt{MATLAB}}
\newcommand{\rmi}{\mathrm{i}\hskip1pt}
\newcommand{\RR}{\mathbb{R}}
\newcommand{\CC}{\mathbb{C}}
\newcommand{\kk}{\tau}

\begin{document}

\title{A $\mu$-mode integrator for solving evolution equations \\ in Kronecker form}
\author[verona]{Marco Caliari} \ead{marco.caliari@univr.it}
\author[trento]{Fabio Cassini\corref{cor1}} \ead{fabio.cassini@unitn.it}
\author[uibk]{Lukas Einkemmer} \ead{lukas.einkemmer@uibk.ac.at}
\author[uibk]{Alexander Ostermann} \ead{alexander.ostermann@uibk.ac.at}
\author[ljll]{Franco Zivcovich} \ead{franco.zivcovich@sorbonne-universite.fr}
\address[verona]{Department of Computer Science, University of Verona, Italy}
\address[trento]{Department of Mathematics, University of Trento, Italy}
\address[uibk]{Department of Mathematics, University of Innsbruck, Austria}
\address[ljll]{Laboratoire Jacques-Louis Lions, Sorbonne University, France}
\cortext[cor1]{Corresponding author}

\begin{abstract}
In this paper, we propose a $\mu$-mode integrator for computing the solution of stiff evolution equations. The integrator is based on a $d$-dimensional splitting approach and uses exact (usually precomputed) one-dimensional matrix exponentials. We show that the action of the exponentials, i.e.~the corresponding batched matrix-vector products, can be implemented efficiently on modern computer systems. We further explain how $\mu$-mode products can be used to compute spectral transforms efficiently even if no \textit{fast} transform is available. We illustrate the performance of the new integrator by solving, among the others, three-dimensional linear and nonlinear Schr\"odinger equations, and we show that the $\mu$-mode integrator  can significantly outperform numerical methods well established in the field. We also discuss how to efficiently implement this integrator on both multi-core CPUs and GPUs. Finally, the numerical experiments show that using GPUs results in performance improvements between a factor of $10$ and $20$, depending on the problem.
\end{abstract}
\maketitle

\section{Introduction}
Due to the importance of simulation in various fields of science and engineering, devising efficient numerical methods for solving evolutionary partial differential equations has received considerable interest in the literature. For linear problems with time-invariant coefficients, after discretizing in space, the task of solving the partial differential equation is equivalent to computing the action of a matrix exponential to a given initial value. Computing the action of matrix exponentials is also a crucial ingredient to devise efficient numerical methods for nonlinear partial differential equations; for example, in the context of exponential integrators~\cite{hochbruck2010} or splitting methods \cite{mclachlan2002}.

Despite the significant advances made in constructing more efficient numerical algorithms, efficiently computing the action of large matrix functions remains a significant challenge. In this paper, we propose a $\mu$-mode integrator that performs this computation for matrices in Kronecker form by computing the action of one-dimensional matrix exponentials only. In $d$ dimensions and with $n$ grid points per dimension the number of arithmetic operations required scales as $\mathcal{O}(n^{d+1})$. Nevertheless, such an approach would not have been viable in the past. With the increasing gap between the amount of floating point operations compared to the amount of memory transactions modern computer systems can perform, however, this is no longer a consequential drawback. In fact, (batched) matrix-matrix multiplications, as are required for this algorithm, can achieve performance close to the theoretical limit of the hardware, and they do not suffer from the irregular memory accesses that plague implementations based on sparse matrix formats. This is particularly true on accelerators, such as Graphic Processing Units (GPUs).
Thus, on modern computer hardware, the proposed method is extremely effective. In this paper, we will show that for a range of problems the proposed $\mu$-mode integrator can outperform well established integrators that are commonly used in the field. We investigate the performances of the method for a two-dimensional pipe flow example. Then, we consider three-dimensional linear Schr\"odinger equations with time dependent and time independent potentials, in combination with Hermite spectral discretization, as well as a cubic nonlinear Schr\"odinger equation (Gross--Pitaevskii equation) in three space dimensions. In this context, we will also provide a discussion on the implementation of the method for multi-core CPUs and GPUs.

The $\mu$-mode integrator is exact for linear problems in \textit{Kronecker form} (see section \ref{sec:mumod} for more details). The discretization of many differential operators with constant coefficients fits into this class (e.g., the Laplacian operator $\Delta$ and the $\rmi \Delta$ operator that is commonly needed in quantum mechanics), as well as some more complicated problems (e.g.~the Hamiltonian for a particle in a harmonic potential). For nonlinear partial differential equations, the approach can be used to solve the part of the problem that is in Kronecker form: for example, in the framework of a splitting method.

The $\mu$-mode integrator is related to dimension splitting schemes such as alternating direction implicit (ADI) schemes (see, e.g., \cite{gasteiger2016,hochbruck2017,namiki1999,peaceman1955}). However, while the main motivation for the dimension splitting in ADI is to obtain one-dimensional matrix equations, for which efficient solvers such as the Thomas algorithm are known, for the $\mu$-mode integrator the main utility of the dimension splitting is the reduction to one-dimensional problems for which matrix exponentials can be computed efficiently. Because of the exactness property described above, for many problems the $\mu$-mode integrator can be employed with a much larger step size compared to implicit methods such as ADI. This is particularly true for highly oscillatory problems, where both implicit and explicit integrators do suffer from small time steps (see, e.g., \cite{ascher1999}).

In the context of spectral decompositions, commonly employed for pseudospectral methods, the structure of the problem also allows us to use $\mu$-mode products to efficiently compute spectral transforms from the space of values to the space of coefficients (and vice versa) even if no $d$-dimensional \textit{fast} transform is available.

The outline of the paper is as follows. In section \ref{sec:mumod} we describe the proposed $\mu$-mode integrator and explain in detail what it means for a differential equation to be in Kronecker form. We also discuss for which class of problems the integrator is particularly efficient. We then show, in section \ref{sec:spectral}, how $\mu$-mode products can be used to efficiently compute arbitrary spectral transforms. Numerical results that highlight the efficiency of the approach will be presented in section \ref{sec:numerical}. The implementation on modern computer architectures, which includes performance results for multi-core CPU and GPU based systems, will be discussed in section \ref{sec:implementation}. Finally, in section \ref{sec:conclusions} we draw some conclusions.

\section{The $\mu$-mode integrator for differential equations in Kronecker form \label{sec:mumod}}

As a simple example that introduces the main idea, we consider the two-dimensional heat equation
\begin{equation}\label{eq:heat2D}
\begin{aligned}
\partial_t u(t,\bx) &= \Delta u(t,\bx) = \left( \partial_1^2 + \partial_2^2\right) u(t,\bx),\qquad \bx\in\Omega\subset\RR^2,\quad t>0,\\
u(0,\bx) &= u_0(\bx)
\end{aligned}
\end{equation}
on a rectangle, subject to appropriate boundary conditions (e.g.~periodic, homogeneous Dirichlet or homogeneous Neumann). Its analytic solution is given by
\begin{equation}\label{eq:heat2D-sol}
u(t,\cdot) = \ee^{t\Delta}u_0 = \ee^{t\partial_1^2}\ee^{t\partial_2^2}u_0 = \ee^{t\partial_2^2}\ee^{t\partial_1^2}u_0,
\end{equation}
where the last two equalities result from the fact that the partial differential operators $\partial_1^2$ and $\partial_2^2$ commute.

Discretizing \eqref{eq:heat2D} by finite differences on a Cartesian grid with $n_1\times n_2$ grid points results in the linear differential equation
\begin{equation}\label{eq:heat2D-fd}
\bu'(t) = \left(I_2\otimes A_1 + A_2 \otimes I_1\right) \bu(t),\quad\bu(0)=\bu_0
\end{equation}
for the unknown vector $\bu(t)$. Here, $A_1$ is a (one-dimensional) stencil matrix for $\partial_1^2$ on the grid points $x_1^{i_1}$, $1\le i_1\le n_1$, and $A_2$ is a (one-dimensional) stencil matrix for $\partial_2^2$ on the grid points $x_2^{i_2}$, $1\le i_2\le n_2$. The symbol $\otimes$ denotes the standard Kronecker product between two matrices.
Since the matrices $I_2\otimes A_1$ and $A_2\otimes I_1$ trivially commute, the solution of \eqref{eq:heat2D-fd} is given by
\begin{equation*}
\bu(t) = \ee^{t\left(I_2\otimes A_1 + A_2 \otimes I_1\right)}\bu_0 = \ee^{tI_2\otimes A_1}\ee^{t A_2 \otimes I_1}\bu_0 = \ee^{t A_2 \otimes I_1}\ee^{tI_2\otimes A_1}\bu_0,
\end{equation*}
which is the discrete analog of \eqref{eq:heat2D-sol}.

Using the tensor structure of the problem, the required actions of the large matrices $\ee^{tI_2\otimes A_1}$ and $\ee^{t A_2 \otimes I_1}$ on a vector can easily be reformulated. Let $\mathbf{U}(t)$ be the order two tensor of
size $n_1\times n_2$ (in fact, a matrix) whose stacked columns form the vector $\bu(t)$. The indices of this matrix reflect the structure of the grid. In particular
\begin{equation*}
\mathbf{U}(t)(i_1,i_2) = u(t,x_1^{i_1}, x_2^{i_2}), \qquad i_1 = 1, \ldots, n_1, \quad i_2 = 1,\ldots, n_2.
\end{equation*}

Using this tensor notation, problem \eqref{eq:heat2D-fd} takes the form
\begin{equation*}
\mathbf{U}'(t) = A_1\, \mathbf{U}(t) + \mathbf{U}(t) A_2^{\sf T},\quad \mathbf{U}(0)= \mathbf{U}_0,
\end{equation*}
and its solution can be expressed as
\begin{equation}\label{eq:sol2D-tensor}
\mathbf{U}(t) =  \ee^{t A_1}\mathbf{U}_0 \,\ee^{t A_2^{\sf T}},
\end{equation}
see \cite{N69}.
From this representation, it is clear that $\mathbf{U}(t)$ can be computed as the action of the small matrices $\ee^{t A_1}$ and $\ee^{t A_2}$ on the tensor $\mathbf{U}_0$. More precisely, the matrices $\ee^{t A_1}$ and $\ee^{t A_2}$ act on the first and second indices of $\mathbf{U}$, respectively. The computation of \eqref{eq:sol2D-tensor} can thus be performed by the simple algorithm
\begin{equation*}
\begin{aligned}
\mathbf{U}^{(0)} &= \mathbf{U}_0,\\
\mathbf{U}^{(1)}(\cdot,i_2) &= \ee^{tA_1}\mathbf{U}^{(0)}(\cdot,i_2),\qquad i_2=1,\ldots,n_2,\\
\mathbf{U}^{(2)}(i_1,\cdot) &= \ee^{tA_2}\mathbf{U}^{(1)}(i_1,\cdot),\qquad i_1=1,\ldots,n_1,\\
\mathbf{U}(t)&=\mathbf{U}^{(2)}.
\end{aligned}
\end{equation*}

It should be duly noted that the $\mu$-mode integrator is not restricted to the simple example considered until now. Indeed, let us consider the differential equation
\begin{equation}\label{eq:discrete}
	\bu'(t) = M \bu(t),\quad \bu(0) = \bu_0,
\end{equation}
where
\begin{equation*}
	M = \sum_{\mu=1}^{d} A_{\otimes \mu}
\end{equation*}
and
\begin{equation}\label{eq:tensor-structure}
	A_{\otimes \mu} = I_d\otimes \cdots\otimes I_{\mu+1}\otimes A_\mu \otimes I_{\mu-1}\otimes \cdots\otimes I_1.
\end{equation}
Here, $A_\mu$ denotes an \textit{arbitrary} $n_\mu\times n_\mu$ matrix while $I_\mu$ is the identity matrix of size $n_\mu$, $1\le \mu\le d$. The matrix $M$ is also known in the literature as the \textit{Kronecker sum} of the matrices $A_\mu$ and is denoted by
\begin{equation*}
	M = A_d \oplus A_{d-1} \oplus \cdots \oplus A_2 \oplus A_1. 
\end{equation*}
Condition \eqref{eq:tensor-structure} holds true for a range of equations with linear and constant coefficient differential operators on tensor product domains. Examples in this class include, after space discretization, the diffusion-advection-absorption equation
	\[ \partial_t u(t,\bx) = \alpha \Delta u(t,\bx) + \beta \cdot \nabla u(t,\bx) - \gamma u(t,\bx) \]
	or the Schr\"odinger equation with potential in Kronecker form
	\[
	\rmi \partial_t \psi(t,\bx) = -\frac{1}{2} \Delta \psi(t,\bx) + \left(\sum_{\mu=1}^d V(t,x_\mu)\right) \psi(t,\bx).
	\]
Condition \eqref{eq:tensor-structure} is fulfilled also for some problems with non-constant coefficient differential operators, see section \ref{sec:pipeflow} for an example. We will consider these and other equations later in the paper to perform numerical examples.

Equation \eqref{eq:discrete} is what we call a linear problem in \emph{Kronecker form}, and its solution is obviously given by
$$
\bu(t) =  \ee^{t A_{\otimes 1}}\cdots \ee^{t A_{\otimes d}} \bu_0,
$$
where the single factors $\ee^{t A_{\otimes \mu}}$
mutually commute.
Again, the computation of $\bu(t)$ just requires the actions of the small matrices $\ee^{t A_\mu}$. More precisely, consider the order $d$
 tensor $\mathbf{U}(t)$ of size $n_1\times \cdots \times n_d$ that collects the values of a function $u$ on a Cartesian grid, i.e.
$$
\mathbf{U}(t)(i_1,\ldots,i_d) = u(t,x_1^{i_1},\ldots,x_d^{i_d}),\qquad 1\le i_\mu\le n_\mu, \quad 1\le \mu \le d.
$$
Then, in the same way as in the two-dimensional heat equation case, the computation of $\bu(t)$ can be performed by
\begin{equation}\label{eq:dd-tensor}
\begin{aligned}
\mathbf{U}^{(0)} &= \mathbf{U}_0,\\[1.5mm]
\mathbf{U}^{(1)}(\cdot,i_2,\ldots,i_d) &= \ee^{tA_1}\mathbf{U}^{(0)}(\cdot,i_2,\ldots,i_d),\qquad 1\le i_\mu\le n_\mu, \quad 2\le \mu \le d,\\
&\cdots\\
\mathbf{U}^{(d)}(i_1,\ldots,i_{d-1},\cdot) &= \ee^{tA_d}\mathbf{U}^{(d-1)}(i_1,\ldots,i_{d-1},\cdot),\qquad 1\le i_\mu\le n_\mu, \quad 1\le \mu \le d-1,\\[1.5mm]
\mathbf{U}(t)&=\mathbf{U}^{(d)}.
\end{aligned}\end{equation}

We remark that scheme \eqref{eq:dd-tensor} can also be useful as a building block for solving nonlinear partial differential equations. In this case, an exponential or splitting scheme would be used to separate the linear part, which is treated exactly by the integrator \eqref{eq:dd-tensor}, from the nonlinear part which is treated in a different fashion. This is useful for a number of problems. For example, when solving the drift-kinetic equations in plasma physics using an exponential integrator \cite{crouseilles2020,crouseilles2018}, Fourier spectral methods are commonly used. While such FFT based schemes are efficient, it is also well known that they can lead to numerical oscillations \cite{filbet2003}.
Using integrator \eqref{eq:dd-tensor} would allow us to choose a more appropriate space discretization while still retaining efficiency. Another example are diffusion-reaction equations with nonlinear reaction terms that are treated using splitting methods (see, e.g., \cite{einkemmer2018bc,einkemmer2015bc,hundsdorfer2013}). In this case scheme \eqref{eq:dd-tensor} would be used to efficiently solve the subflow corresponding to the linear diffusion.
We further note that a related approach was pursued by \cite{NWZL08} in order to produce schemes that solve two- and three-dimensional biological models.

Implementing integrator \eqref{eq:dd-tensor} requires the computation of $d$ small exponentials of sizes $n_1 \times n_1$, \ldots, $n_d \times n_d$, respectively. If a marching scheme with \emph{constant} time step is applied to (\ref{eq:discrete}), then these matrices can be precomputed once and for all, and their storage cost is negligible compared to that required by the solution $\mathbf U(t)$. Otherwise, we need to compute at every time step new matrix exponentials, whose computational cost still represents only a small fraction of the entire algorithm (see section \ref{sec:heat}). Indeed, the main component of the final cost is represented by the computation of matrix-matrix products of size $n_\mu \times n_\mu$ times $n_\mu \times (n_1  \cdots n_{\mu-1} n_{\mu+1}  \cdots n_{d})$. Thus, the computational complexity of the algorithm is $\mathcal{O}(N\max_\mu n_\mu)$, where $N=n_1 \cdots n_d$ is the total number of degrees of freedom.

Clearly, we can solve equation \eqref{eq:discrete} also by directly computing the vector $\ee^{t M}\mathbf{u}_0$. In fact $M$ is an $N\times N$ sparse matrix and, when it is too large for the explicit computation of $\ee^{t M}$, the action of the matrix exponential can be approximated by polynomial methods such as Krylov projection (see, for instance, \cite{GRT18,NW12}), Taylor series \cite{AMH11}, or polynomial interpolation (see, for instance, \cite{CCZ20, CKOR16,CKZ18}). All these iterative methods require one matrix-vector product per iteration, which costs $\mathcal{O}(N)$ plus additional vector operations. The number of iterations, however, highly depends on the norm and some properties  of the matrix, such as the normality, the condition number, and the stiffness, and it is not easy to predict it. Moreover, for Krylov methods, one has to take into account the storage of a full matrix with $N$ rows and as many columns as the dimension of the Krylov subspace.

Also, an implicit scheme based on a Krylov solver could be applied to integrate equation~\eqref{eq:discrete}. In particular, if we restrict our attention to the heat equation case and the conjugate gradient method, for example, $\mathcal{O}(\max_\mu n_\mu)$ iterations are needed for the solution (see the convergence analysis in \cite[Chap. 6.11]{saad2003}), and each iteration requires a sparse matrix-vector product which is $\mathcal{O}(N)$. Hence, the resulting computational complexity is the same as for the proposed algorithm.  However, on modern hardware architectures memory transactions are much more costly than performing floating point operations. A modern CPU or GPU can easily perform many tens of arithmetic operations in the same time it takes to read/write a single number from/to memory (see the discussion in section \ref{sec:implementation}).

Summarizing, our scheme has the following advantages:
\begin{itemize}
\item For a heat equation the proposed integrator only requires $\mathcal{O}(N)$ memory operations, compared to an
implicit scheme which requires $\mathcal{O}(N\max_\mu n_\mu)$ memory operations. This has huge performance implications on all modern computer architectures. For other classes of PDEs the analysis is more complicated. However, in many situations similar results can be obtained.
    \item Very efficient implementations of matrix-matrix products that operate close to the limit of the hardware are available. This is not the case for iterative schemes which are based on sparse matrix-vector products.
    \item The computation of pure matrix exponentials of small matrices is less prone to the problems that affect the approximation of the action of the
     (large) matrix exponential.
    \item The proposed integrator is often able to take much larger time step sizes than, for example, an ADI scheme, as it computes the exact result for equations in Kronecker form.
    \item Conserved quantities of the underlying system, such as mass, are preserved by the integrator.
\end{itemize}
    We will in fact see that the proposed integrator can outperform algorithms with linear computational complexity (see sections \ref{sec:schrti} and \ref{sec:schrtd}).

Equation \eqref{eq:dd-tensor} gives perhaps the most intuitive picture of the proposed approach. However, we can also formulate this problem in terms of $\mu$-fibers. Indeed, let $\mathbf{U}\in\CC^{n_1\times\cdots\times n_d}$ be an order $d$ tensor. A \emph{$\mu$-fiber} of $\mathbf{U}$ is a vector in $\CC^{n_\mu}$ obtained by fixing every index of the tensor but the $\mu$th. In these terms,  $\mathbf{U}^{(\mu-1)}(i_1,\ldots,i_{\mu-1},\cdot,i_{\mu+1},\ldots,i_d)$ is a $\mu$-fiber of the tensor $\mathbf{U}^{(\mu-1)}$, and every line in formula (\ref{eq:dd-tensor}) corresponds to the action of the matrix $\ee^{tA_\mu}$ on the $\mu$-fibers of $\mathbf{U}^{(\mu-1)}$. By means of $\mu$-fibers, it is possible to define the following operation.

\begin{definition} \label{def:mma}
Let $L \in \CC^{m\times n_\mu}$ be a matrix. Then the \emph{$\mu$-mode product}\footnote{Also known as mode-$n$ product, $n$-mode product or mode-$\alpha$ multiplication, depending on the convention.} of $L$ with $\mathbf{U}$, denoted by $\mathbf S = \mathbf{U} \times_\mu L$, is the tensor $\mathbf S\in\CC^{n_1\times\ldots\times n_{\mu-1}\times m\times n_{\mu+1}\times\ldots \times n_d}$ obtained by multiplying the matrix $L$ onto the $\mu$-fibers of $\mathbf{U}$, that is
\begin{equation*}
\mathbf{S}(i_1, \cdots, i_{\mu-1}, i, i_{\mu+1}, \cdots, i_d)= \sum_{j=1}^{n_\mu} L_{ij} \mathbf{U}(i_1, \cdots, i_{\mu-1},j, i_{\mu+1}, \cdots, i_d), \qquad 1\le i\le m.
\end{equation*}
\end{definition}\vspace{-5mm}

According to this definition, it is clear that in formula (\ref{eq:dd-tensor}) we are performing $d$ consecutive $\mu$-mode products with the matrices $\ee^{tA_\mu}$, $1\le \mu \le d$. We can therefore write scheme \eqref{eq:dd-tensor} as follows
\[ \mathbf{U}(t) = \mathbf{U}_0\times_1 \ee^{tA_1} \times_2\ldots\times_d \ee^{tA_d}. \]
This is the reason why we call the proposed method the $\mu$-mode integrator. Notice that the concatenation of $\mu$-mode products of $d$ matrices with a tensor is also known as the \emph{Tucker operator} (see \cite{K06}), and it can be performed using efficient level-3 BLAS operations. For more information on tensor algebra and the $\mu$-mode product we refer the reader to \cite{KB09}.

\section{Application of the $\mu$-mode product to spectral decomposition and reconstruction \label{sec:spectral}}

Problems of quantum mechanics with vanishing boundary conditions are often set in an unbounded spatial domain. In this case, the spectral decomposition in space by Hermite functions is appealing (see~\cite{BS05,TCN09}), since it allows to treat boundary conditions  in a natural way (without imposing artificial periodic boundary conditions as required by Fourier spectral methods, for example).

Consider the multi-index $\mathbf i=(i_1,\ldots,i_d)\in\mathbb{N}_0^d$ and the coordinate vector $\mathbf x=(x_1,\ldots,x_d)$ belonging to $\RR^d$. We define the $d$-variate functions
$\mathcal{H}_{\mathbf{i}}(\mathbf x)$ as
\begin{equation*}
    \mathcal{H}_{\mathbf{i}}(\mathbf x)=\prod_{\mu=1}^d H_{i_\mu}(x_\mu)\ee^{-x_\mu^2/2},
\end{equation*}
where $\{H_{i_\mu}(x_\mu)\}_{i_\mu}$ is the family of Hermite polynomials ortho\emph{normal} with respect to the weight function $\ee^{-x_\mu^2}$ on $\RR$, that is
\begin{equation*}
    \int_{\RR^d} \mathcal{H}_{\mathbf i}(\mathbf x)
    \mathcal{H}_{\mathbf{j}}(\mathbf x)
    d \mathbf x=
    \delta_{\mathbf {i}\mathbf {j}}.
\end{equation*}
We recall that Hermite functions satisfy
\begin{equation*}
\left(-\frac{1}{2}\sum_{\mu=1}^d(\partial_{\mu}^2-x_\mu^2)\right)
\mathcal{H}_{\mathbf i}(\mathbf x)=\lambda_{\mathbf i}
\mathcal{H}_{\mathbf i}(\mathbf x),
\end{equation*}
where
\begin{equation*}
\lambda_{\mathbf i}=\sum_{\mu=1}^d\left(\frac{1}{2}+i_\mu\right).
\end{equation*}

In general, we can consider a family of functions
$\phi_{\mathbf i}\colon R_1\times \cdots \times R_d\to \CC$ in tensor form
\begin{equation*}
\phi_{\mathbf i}(\mathbf x) = \prod_{\mu=1}^d \phi_{i_\mu}^{\mu}(x_\mu)
\end{equation*}
which are orthonormal on the Cartesian product of intervals $R_1,\ldots, R_d$ of $\RR$.

If a function $f$ can be expanded into a series
\begin{equation*}
    f(\mathbf x)=\sum_{\mathbf i}
    f_{\mathbf i}\phi_{\mathbf i}(\mathbf x),\qquad
    f_{\mathbf i}\in\CC,
\end{equation*}
then its $\mathbf i$th coefficient is
\begin{equation*}
    f_{\mathbf i}=\int_{R_1\times\cdots\times R_d}f(\mathbf x)
    \overline{\phi_{\mathbf i}}(\mathbf x)d\mathbf x.
\end{equation*}
In order to approximate the integral on the right-hand side, we rely on a tensor-product quadrature formula.
To do so, we consider for each direction $\mu$ a set of $m_\mu$ uni-variate quadrature nodes $X_{\ell_\mu}^\mu$ and weights $W_{\ell_\mu}^\mu$, $0\le \ell_\mu\le m_\mu$, and fix to $k_\mu$ the number of uni-variate functions $\phi_{i_\mu}^{\mu}(x_\mu)$ to be considered.
We have then 
\begin{equation}\label{eq:quadrature}
    \hat f_{\mathbf i}=
    \sum_{\boldsymbol \ell < \mathbf m}f(\mathbf x_{\boldsymbol{\ell}})
    \overline{\phi_{\mathbf i}}(\mathbf x_{\boldsymbol \ell})w_{\boldsymbol \ell},
    \qquad \mathbf i < \mathbf k,
\end{equation}
where $\mathbf x_{\boldsymbol \ell}=(X_{\ell_1}^1,\cdots,X_{\ell_d}^d) \in \RR^d$,
$w_{\boldsymbol \ell}=\prod_{\mu=1}^d W_{\ell_\mu}^\mu$ and
$\mathbf k$ is the multi-index which collects the values $\{k_\mu\}_\mu$.
We show now how $\mu$-mode products can be employed to compute the coefficients of the spectral decomposition
\begin{equation}\label{eq:truncseries}
  \hat f(\mathbf x)=\sum_{\mathbf i<\mathbf k}\hat f_{\mathbf i}
  \phi_{\mathbf i}(\mathbf x)\approx f(\mathbf x)
\end{equation}
of a $d$-variate function and its evaluation on a Cartesian grid.
First of all, for each fixed $\mu$, $1\le \mu \le d$, we define the matrix $\Phi_\mu\in \CC^{k_\mu \times m_\mu}$ with components
\begin{equation*} \label{eq: r2smat}
    (\Phi_\mu)_{i\ell} = \overline{\phi^\mu_{i}}(X_{\ell}^\mu),
\end{equation*}
and we denote by $\mathbf F_{\mathbf{W}}\in\CC^{m_1\times\cdots\times m_d}$ the tensor with elements
$f(\mathbf x_{\boldsymbol \ell})w_{\boldsymbol{\ell}}$  and by $\hat {\mathbf F}\in\CC^{k_1\times\cdots\times k_d}$ the tensor with elements
$\hat f_{\mathbf i}$. Then, in terms of the Tucker operator, we can write equation \eqref{eq:quadrature} as follows
\begin{equation}\label{eq:real2spec}
\hat {\mathbf F}=
\mathbf F_{\mathbf W}\times_1 \Phi_1
\times_2\cdots\times_d \Phi_d.
\end{equation}

It is then possible to evaluate the function $\hat f(\mathbf{x})$
in \eqref{eq:truncseries} at a Cartesian grid $\mathbf{y_p} = (Y^1_{p_1}, \ldots, Y^d_{p_d})$, that is
\begin{equation} \label{eq:collocation}
    \hat f(\mathbf{y_p}) =  \sum_{\mathbf i < \mathbf k}\hat f_{\mathbf i}\phi_{\mathbf i}(\mathbf{y_p}), \qquad \mathbf p < \mathbf q,
\end{equation}
by the Tucker operator, too. Here the component $q_\mu$ of the multi-index $\mathbf{q}$ is the number of uni-variate evaluation points $Y^\mu_{p_\mu}$.
Indeed, if we collect the elements $\hat f(\mathbf y_{\mathbf p})$ in the tensor $\hat{\hat {\mathbf F}} \in \CC^{q_1 \times \cdots \times q_d}$ and, for fixed $\mu$, we define the matrix $\Psi_\mu\in \CC^{q_\mu \times k_\mu}$ with components $(\Psi_\mu)_{pi} = \phi^\mu_{i}(Y_{p}^\mu)$, then 
\begin{equation}\label{eq:spec2real}
    \hat{\hat {\mathbf F}} = \hat {\mathbf F}\times_1 \Psi_1 \times_2 \cdots \times_d\Psi_d
\end{equation}
is the tensor formulation of formula (\ref{eq:collocation}).

Now, we restrict our attention to the common case where the quadrature nodes are chosen in such a way that
\begin{equation*}
\sum_{\boldsymbol \ell < \mathbf m}\phi_{\mathbf i}(\mathbf x_{\boldsymbol \ell})
\overline{\phi_{\mathbf j}}(\mathbf x_{\boldsymbol \ell})w_{\boldsymbol \ell} = \delta_{\mathbf i \mathbf j}, \qquad \mathbf i, \mathbf j < \mathbf k
\end{equation*}
with $\mathbf m=\mathbf k$, that is, the orthonormality relation among the $\phi_{\mathbf i}$ functions is true also at the discrete level.
This is the case, for instance, when using
Gauss--Hermite quadrature nodes for
$\phi_{\mathbf i}(\mathbf x)=\mathcal{H}_{\mathbf i}(\mathbf x)$.
Then, the matrices $\Phi_{\mu}\in\CC^{m_\mu \times m_\mu}$ turn out to be square
and formula \eqref{eq:real2spec} is the \emph{spectral transform} from
the space of values to the space of coefficients.
Moreover, if the evaluation points coincide with the quadrature nodes, then we have $\Psi_{\mu} = \Phi_\mu^*$, where the symbol $*$ denotes the conjugate transpose of the matrix, and formula \eqref{eq:spec2real} is the
\emph{inverse spectral transform} from the space of coefficients
to the space of values.

As mentioned at the beginning of the section, we will employ the Hermite spectral decomposition in some of our experiments (see sections \ref{sec:schrti} and \ref{sec:schrtd}). Hence, we will use (\ref{eq:real2spec}) and (\ref{eq:spec2real}) for the required spectral transforms.

We also remark that a similar approach was pursued in \cite{hashemi2017chebfun} in the framework of three-dimensional Chebyshev interpolation.

\section{Numerical comparison \label{sec:numerical}}

In this section, we will compare the proposed $\mu$-mode integrator with some widely used techniques to solve partial differential equations. For that purpose a range of PDEs, mainly from quantum mechanics, is considered. Concerning the experiments in sections \ref{sec:heat}, \ref{sec:pipeflow} and \ref{sec: gpe}, we will test the proposed method against the following iterative schemes commonly employed to compute the action of the matrix exponential $\ee^{tM}$:
\begin{itemize}
    \item \texttt{expmv}: a polynomial method described in \cite{AMH11} which is based on a Taylor expansion of the exponential;
    \item \texttt{phipm}: a full Krylov method presented in \cite{NW12};
    \item \texttt{kiops}: a Krylov method based on an incomplete orthogonalization process, described in \cite{GRT18}.
\end{itemize}
The \matlab{} source code of these methods is publicly available. Although the underlying algorithms of these schemes only require the action of the matrix on a vector, only \texttt{kiops} is readily available to do that. Therefore, in order to ensure a fair comparison, we feed the functions with the matrix. Moreover, considering the action of the matrix on a vector (which in our case could be performed entirely in tensor formulation by means of sums of $\mu$-mode products) instead of the matrix itself would not result in a speedup for the schemes (see section \ref{sec:heat}).
The tolerance for all the algorithms considered has been set to $2^{-53}$, which corresponds to the machine epsilon for double precision computations. As a measure of cost, we consider the computational time (wall-clock time) needed to solve numerically the differential equation under consideration up to a fixed final time. As mentioned in section \ref{sec:mumod}, the $\mu$-mode integrator requires the explicit computation of small matrix exponentials. This is performed using the internal \matlab\textsuperscript{\textregistered} function \texttt{expm}, which is based on the scaling and squaring rational Pad\'e approximation described in \cite{AMH09}. In this context, another method which could be directly used in \matlab{} is \texttt{exptayotf} from \cite{CZ19}. It is based on a backward stable Taylor approximation for the matrix exponential and is faster than \texttt{expm}. Moreover, as it works in single, double and variable precision arithmetic data types, it produces approximations with the desired accuracy. This is not possible for the iterative schemes which approximate the action of $\ee^{tM}$, because the \matlab\textsuperscript{\textregistered} sparse format is restricted to double precision. Another fast method using a similar technique and suited for double precision is \texttt{expmpol} from \cite{SID18}. We will demonstrate that our \matlab{} implementation of the proposed $\mu$-mode integrator outperforms all the other schemes by at least a factor of 7.

Concerning the experiments in sections \ref{sec:schrti} and  \ref{sec:schrtd}, we compare our $\mu$-mode based approach with a splitting scheme/FFT based space discretization that is well established and efficient. In order to perform direct and inverse Fourier transforms, we employ the internal \matlab\textsuperscript{\textregistered} functions \texttt{fftn} and \texttt{ifftn} respectively, which are in turn based on the very efficient FFTW library \cite{frigo1998fftw}. Care has been taken to ensure that comparisons conducted in \texttt{MATLAB}\textsuperscript{\textregistered} give a good indication of the performance that would be obtained in a compiled language. This is possible here as the majority part of the computational time is spent in the FFT routines. For these problems, we will show that the $\mu$-mode integrator can reach a speedup of at least 5.

All the tests in this section have been conducted on an Intel Core i7-5500U CPU with 12GB of RAM using \matlab\textsuperscript{\textregistered} \texttt{R2020b}.

\subsection{Code validation \label{sec:heat}}

As an introductory test problem, in order to highlight some qualities of our $\mu$-mode method, we consider the three-dimensional heat equation
\begin{equation} \label{eq:he}
\left\{
\begin{aligned}
\partial_t u(t,\mathbf{x})& = \Delta u(t,\mathbf{x}),\qquad \mathbf{x}\in [0,2\pi)^3,\quad t\in[0,T], \\
u(0,\mathbf{x}) &= \cos x_1 + \cos x_2 + \cos x_3
\end{aligned}
\right.
\end{equation}
with periodic boundary conditions.

The equation is discretized in space using centered finite differences with $n_\mu$ grid points in the $\mu$th direction (the total number of degrees of freedom stored in computer memory is hence equal to $N = n_1n_2n_3$). By doing so we obtain the following ordinary differential equation (ODE)
\begin{equation}\label{eq: lineareq}
\mathbf{u}'(t) = M\mathbf{u}(t),
\end{equation}
where $\mathbf{u}$ denotes the vector in which the degrees of freedom are assembled. The exact solution of equation (\ref{eq: lineareq}) is given by the action of the matrix exponential
\begin{equation}\label{eq:exactsol}
\mathbf{u}(t) = \ee^{tM}\mathbf{u}(0).
\end{equation}
The matrix $M$ has the following Kronecker structure
\begin{equation*}
M = I_3 \otimes I_2 \otimes A_1 + I_3 \otimes A_2 \otimes I_1 + A_3 \otimes I_2 \otimes I_1,
\end{equation*}
where $A_\mu \in \RR^{n_\mu\times n_\mu}$ results from the one-dimensional discretization of the operator $\partial^2_{\mu}$, and $I_\mu\in\RR^{n_\mu\times n_\mu}$ is the identity matrix.
The quantity $\mathbf{u}(t)$ can be seen as vectorization of the tensor $\mathbf{U}(t)$, and we can write (\ref{eq:exactsol}) in tensor form as
\begin{equation*}
\mathbf{U}(t) = \mathbf{U}(0) \times_1 \ee^{tA_1} \times_2 \ee^{tA_2} \times_3 \ee^{tA_3},
\end{equation*}
where $\mathbf U(t)(i_1,i_2,i_3)=\mathbf u(t)_{i_1+n_1(i_2-1)+n_1n_2(i_3-1)}.$

We now present three numerical tests.
\begin{description}
    \item[Test 1.] We consider second order centered finite differences and compute the solution at time $T=1$ for $n_{\mu} = n$, $\mu=1,2,3$ with various $n$.  We investigate the wall-clock time as a function of the problem size.
    \item[Test 2.] We fix the problem size ($n_\mu=40$, $\mu = 1,2,3$) and compute the solution at time $T=1$ for different orders $p$ of the finite difference scheme. We thereby investigate the wall-clock time as a function of the sparsity pattern of $M$.
    \item[Test 3.]  We consider second order centered finite differences and fix the problem size ($n_\mu=40$, $\mu = 1,2,3$). We then compute the solution at different final times $T$. By doing so we investigate the wall-clock time as a function of the norm of $M$.
\end{description}

\begin{figure}
    \centering
%
%
\definecolor{mycolor1}{rgb}{0.00000,0.44700,0.74100}%
\definecolor{mycolor2}{rgb}{0.85000,0.32500,0.09800}%
\definecolor{mycolor3}{rgb}{0.92900,0.69400,0.12500}%
\definecolor{mycolor4}{rgb}{0.49400,0.18400,0.55600}%
\begin{tikzpicture}

\begin{axis}[%
small,
width=1.232in,
height=3.0946in,
at={(0.78in,0.492in)},
scale only axis,
xmin=1,
xmax=5,
xtick={1,2,3,4,5},
xticklabels={{40},{55},{70},{85},{100}},
xlabel style={font=\color{white!15!black}},
xlabel style={font=\footnotesize},
xlabel={$n$},
ymode=log,
ymin=0.0001,
ymax=200,
yminorticks=true,
ylabel style={font=\color{white!15!black}},
ylabel style={font=\footnotesize},
ylabel={Wall-clock time (s)},
axis background/.style={fill=white},
]
\addplot [thick, color=mycolor1, mark=x, mark options={solid, mycolor1, scale = 1.5}, forget plot]
  table[row sep=crcr]{%
1	1.06736625\\
2	6.6345089\\
3	22.7806582\\
4	57.7853785\\
5	112.631793\\
};
\addplot [thick, color=mycolor2, mark=+, mark options={solid, mycolor2, scale = 1.5}, forget plot]
  table[row sep=crcr]{%
1	0.2817965\\
2	1.43677\\
3	5.5259876\\
4	12.64119\\
5	29.773889\\
};
\addplot [thick, color=mycolor3, mark=asterisk, mark options={solid, mycolor3, scale = 1.5}, forget plot]
  table[row sep=crcr]{%
1	0.1291513\\
2	0.563185\\
3	1.6472958\\
4	6.0901885\\
5	9.599881\\
};
\addplot [thick, color=mycolor4, mark=o, mark options={solid, mycolor4, scale = 1.5}, forget plot]
  table[row sep=crcr]{%
1	0.0013069955\\
2	0.002412247\\
3	0.007108788\\
4	0.014331325\\
5	0.02008297\\
};
\end{axis}

\begin{axis}[%
small,
width=1.232in,
height=3.0946in,
at={(2.465in,0.492in)},
scale only axis,
xmin=2,
xmax=10,
xtick={2,4,6,8,10},
xticklabels={{2},{4},{6},{8},{$\infty$}},
xlabel style={font=\color{white!15!black}},
xlabel style={font=\footnotesize},
xlabel={$p$},
ymode=log,
ymin=0.0001,
ymax=200,
yminorticks=true,
axis background/.style={fill=white},
]
\addplot [thick, color=mycolor1, mark=x, mark options={solid, mycolor1, scale = 1.5}, forget plot]
  table[row sep=crcr]{%
2	1.05226185\\
4	2.42796925\\
6	3.5979371\\
8	5.6869289\\
10	39.130831\\
};
\addplot [thick, color=mycolor2, mark=+, mark options={solid, mycolor2, scale = 1.5}, forget plot]
  table[row sep=crcr]{%
2	0.2663221\\
4	0.51102825\\
6	0.7438609\\
8	0.97310315\\
10	7.809532\\
};
\addplot [thick, color=mycolor3, mark=asterisk, mark options={solid, mycolor3, scale = 1.5}, forget plot]
  table[row sep=crcr]{%
2	0.12988305\\
4	0.2016072\\
6	0.2682686\\
8	0.2934256\\
10	1.717691\\
};
\addplot [thick, color=mycolor4, mark=o, mark options={solid, mycolor4, scale = 1.5}, forget plot]
  table[row sep=crcr]{%
2	0.001347043\\
4	0.0011753645\\
6	0.0015866315\\
8	0.001137694\\
10	0.00142465\\
};
\end{axis}

\begin{axis}[%
small,
width=1.232in,
height=3.0946in,
at={(4.15in,0.492in)},
scale only axis,
xmin=-2,
xmax=2,
xtick={-2,-1,0,1,2},
xticklabels={{$\text{2}^{\text{-2}}$},{$\text{2}^{\text{-1}}$},{$\text{2}^{\text{0}}$},{$\text{2}^{\text{1}}$},{$\text{2}^{\text{2}}$}},
xlabel style={font=\color{white!15!black}},
xlabel style={font=\footnotesize},
xlabel={$T$},
ymode=log,
ymin=0.0001,
ymax=200,
yminorticks=true,
axis background/.style={fill=white},
legend style={legend cell align=left, align=left, draw=white!15!black},
legend style={font=\footnotesize}
]
\addplot [thick, color=mycolor1, mark=x, mark options={solid, mycolor1, scale = 1.5}]
  table[row sep=crcr]{%
-2	0.2690752\\
-1	0.5695619\\
0	1.0381528\\
1	2.02615065\\
2	3.93648955\\
};
\addlegendentry{expmv}

\addplot [thick, color=mycolor2, mark=+, mark options={solid, mycolor2, scale = 1.5}]
  table[row sep=crcr]{%
-2	0.09315235\\
-1	0.1277262\\
0	0.24774025\\
1	0.474892\\
2	0.8246948\\
};
\addlegendentry{phipm}

\addplot [thick, color=mycolor3, mark=asterisk, mark options={solid, mycolor3, scale = 1.5}]
  table[row sep=crcr]{%
-2	0.0564548\\
-1	0.0881799\\
0	0.12201685\\
1	0.1818998\\
2	0.2460607\\
};
\addlegendentry{kiops}

\addplot [thick, color=mycolor4, mark=o, mark options={solid, mycolor4, scale = 1.5}]
  table[row sep=crcr]{%
-2	0.001267674\\
-1	0.001173257\\
0	0.0011590785\\
1	0.0011509815\\
2	0.001157699\\
};
\addlegendentry{$\mu\text{-mode}$}

\end{axis}
\end{tikzpicture}%
    \caption{The wall-clock time for solving the heat equation \eqref{eq:he} is shown as a function of $n$ (left), of the order of the finite difference scheme $p$ (middle), and of the final time $T$ (right). Note that $p=\infty$ corresponds to a spectral space discretization.} \label{fig: expmatvec}
\end{figure}
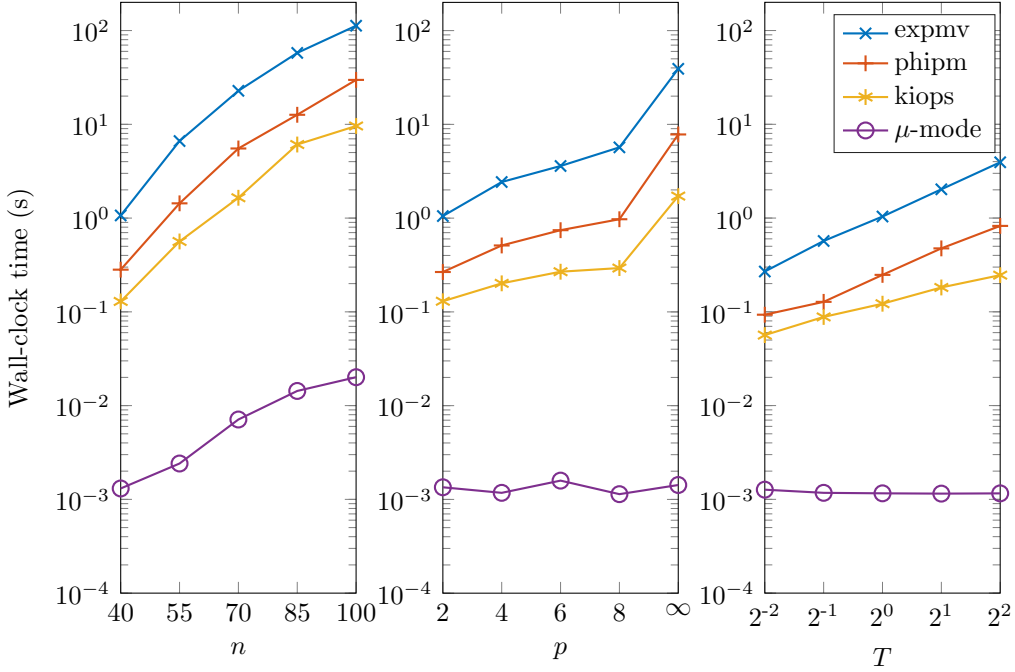

The corresponding results are shown in Figure \ref{fig: expmatvec}. We see that the proposed $\mu$-mode integrator is always the fastest algorithm. The difference in computational time is at least a factor of $60$.

Concerning the first test, we measure also the relative error between the analytical solution and the numerical one. As the dimensional splitting performed by the $\mu$-mode integrator is exact, its errors are equal to the ones obtained by computing (\ref{eq:exactsol}) using the other algorithms. Indeed, for the values of $n$ under consideration, we obtain 2.06e-03, 1.09e-03, 6.71e-04, 4.55e-04 and 3.29e-04 for all the methods. We highlight also that the main cost of the $\mu$-mode integrator is represented by the computation of the $\mu$-mode products and not by the exponentiation of the matrices $A_\mu$  (see Table \ref{tab: mmpbd}). Lastly, notice that the iterative algorithms would not have taken advantage from the computation of the internal matrix-vector products, which constitute their main cost, in tensor formulation (i.e. by means of sums of $\mu$-mode products). Indeed, if we measure the wall-clock time for a single action of the matrix on a vector we observe, for the values of $n$ under consideration, a speedup of averagely 1.5 times by using the standard sparse matrix-vector product as opposed to the tensor formulation.

\begin{table}[!ht]
    \centering
    \bgroup\def\arraystretch{1.3}
    \begin{tabular}{c||ccccc}
        $n$ & $40$ & $55$ & $70$ & $85$ & $100$ \\
        \hline
        \texttt{expm}&  0.52  & 0.71 & 1.37 &  3.15 & 3.54 \\
        $\mu$-mode products& 0.79  & 1.71 &  5.74  & 10.92 & 16.89 \\
        \hline
        Total & 1.31  & 2.42  & 7.11 &  14.07 & 20.43
    \end{tabular}
    \egroup
    \caption{Breakdown of wall-clock time (in ms) for the $\mu$-mode integrator for different values of $n$ (cf. left plot of Figure \ref{fig: expmatvec}).}\label{tab: mmpbd}
\end{table}

The second test shows that the iterative schemes see a decrease in performance when decreasing the sparsity of the matrix (i.e.~by increasing the order of the method $p$ or by using a spectral approximation). This effect is particularly visible when performing a spectral discretization, which results in full matrices  $A_\mu$. On the other hand, the $\mu$-mode integrator is largely unaffected as it computes the exponential of the \textit{full} matrices $A_\mu$, independently of the initial sparsity pattern, by using \texttt{expm}.

Similar observations can be made for the third test. While the iterative schemes suffer from increasing computational time as the norm of the matrix increases, for the $\mu$-mode integrator this is not the case. The reason for this is that the scaling and squaring algorithm in \texttt{expm} scales very favorably as the norm of the matrix increases.

\subsection{Pipe flow \label{sec:pipeflow}}
To demonstrate that the $\mu$-mode integrator can be used for some problems with non-constant coefficients, we consider a model for a fluid flowing in a pipe. The main assumptions are that of radial symmetry (i.e.~the solution does not depend on the angle variable in the circular cross section, see for example \cite{taylor1953dispersion}) and a prescribed length-dependent flow velocity. In this case we obtain the following diffusion-advection equation for the concentration $c$
\begin{equation} \label{eq:pf}
        \partial_t c(t,\rho,z) = \alpha\left( \partial_{\rho\rho}c(t,\rho,z) + \frac{1}{\rho}\partial_{\rho}c(t,\rho,z) +\partial_{zz}c(t,\rho,z)\right) - s(z)\partial_{z}c(t,\rho,z),
\end{equation}
where $t\in[0,T]$, $\rho\in[\rho_{\min},\rho_{\max}]$ and $z\in[0,z_{\max}]$. Here $\alpha$ is the diffusivity and $s(z)$ represents the advection velocity.

After space discretization, which in our case is performed by means of second order centered finite differences with equal number of discretization points $n_\mu$ in each direction (i.e.~$n_\mu=n$, with $\mu=1,2$), the resulting ODE is a linear problem in Kronecker form \eqref{eq:tensor-structure}. The system can then be integrated exactly by the $\mu$-mode integrator.

For the simulations conducted, we use the following initial and boundary conditions
\begin{equation*}
	\left\{
	\begin{aligned}
		& c(0,\rho,z) = \exp(-8(\rho-\rho_0)^2 - 8(z-z_0)^2), \\
		& c(t,\rho,0) = 0, \\
		&\partial_zc(t,\rho,z_{\max}) = 0, \\
		&\partial_\rho c(t,\rho_{\min},z) = 0, \\
		&\partial_\rho c(t,\rho_{\max},z) = 0, \\
		\end{aligned}
	\right.
\end{equation*}
 while the flow velocity is set to
\begin{equation*}
	s(z) = 2 + \tanh(4(z-5/2))-\tanh(4(z-5)).
\end{equation*}
The parameters are chosen as $\rho_{\min} = 0.1$, $\rho_{\max} = 5$, $z_{\max} = 8$, $\alpha = 1/90$, $\rho_0 = (\rho_{\min} + \rho_{\max})/2$ and $z_0 = 3/2$. The structure of the problem does not allow an effective use of FFT based methods.

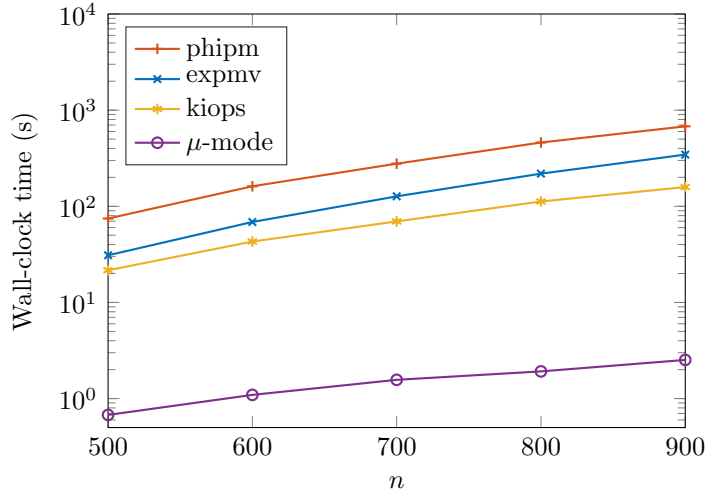
\begin{figure}
	\centering
%
%
\definecolor{mycolor1}{rgb}{0.92941,0.69412,0.12549}%
\definecolor{mycolor2}{rgb}{0.85098,0.32549,0.09804}%
\definecolor{mycolor3}{rgb}{0.00000,0.44706,0.74118}%
\definecolor{mycolor4}{rgb}{0.49412,0.18431,0.55686}%
\begin{tikzpicture}

\begin{axis}[%
width=3.021in,
height=2.166in,
at={(0.758in,0.481in)},
scale only axis,
xmin=1,
xmax=5,
xtick={1,2,3,4,5},
xticklabels={{500},{600},{700},{800},{900}},
xlabel style={font=\color{white!15!black}},
xlabel style={font=\footnotesize},
xticklabel style={font=\footnotesize},
xlabel={$n$},
ymode=log,
ymin=0.5,
ymax=10000,
yminorticks=true,
ylabel style={font=\color{white!15!black}},
ylabel={Wall-clock time (s)},
ylabel style={font=\footnotesize},
yticklabel style={font=\footnotesize},
axis background/.style={fill=white},
legend style={at={(0.03,0.97)}, anchor=north west, legend cell align=left, align=left, draw=white!15!black},
legend style={font=\footnotesize}
]
\addplot [thick, color=mycolor2, mark=+, mark options={solid, mycolor2}]
table[row sep=crcr]{%
	1	74.5962\\
	2	161.4983\\
	3	277.6028\\
	4	460.6930\\
	5	678.6271\\
};
\addlegendentry{phipm}

\addplot [thick, color=mycolor3, mark=x, mark options={solid, mycolor3}]
  table[row sep=crcr]{%
	1	30.9257\\
2	68.6897\\
3	126.8846\\
4	218.5680\\
5	344.7393\\
};
\addlegendentry{expmv}

\addplot [thick, color=mycolor1, mark=asterisk, mark options={solid, mycolor1}]
  table[row sep=crcr]{%
1	21.6231\\
2	42.9917\\
3	69.6489\\
4	112.2069\\
5	158.4907\\
};
\addlegendentry{kiops}

\addplot [thick, color=mycolor4, mark=o, mark options={solid, mycolor4}]
  table[row sep=crcr]{%
1	0.6763\\
2	1.0928\\
3	1.5674\\
4	1.9183\\
5	2.5244\\
};
\addlegendentry{$\mu\text{-mode}$}

\end{axis}
\end{tikzpicture}%
	\caption{Wall-clock time (in seconds) for the integration of (\ref{eq:pf}) up to $T=4$ as a function of $n$ (total number of degrees of freedom $N=n^2$).}
	\label{fig: pipeflow}
\end{figure}

The results of the experiment are presented in Figure \ref{fig: pipeflow}. The $\mu$-mode integrator outperforms all the iterative methods by a consistent factor, with an average speedup of 45 times with respect to \texttt{kiops}, the fastest competitor in this simulation. 

\subsection{Schr\"odinger equation with time independent potential \label{sec:schrti}}

In this section we solve the Schr\"odinger equation in three space dimensions
\begin{equation}\label{eq: schtind}
\left\{\begin{aligned}
\rmi\partial_t\psi(t,\mathbf x)&=
-\frac{1}{2}\Delta\psi(t,\mathbf x)
+V(\mathbf x)\psi(t,\mathbf x), \qquad \mathbf x\in\RR^3, \quad t\in[0,1] \\
\psi(0,\mathbf x)&=\psi_0(\mathbf x)
\end{aligned}\right.
\end{equation}
with a time independent potential $V(\mathbf x)=V_1(x_1)+V_2(x_2)+V_3(x_3)$, where
\begin{equation*}
V_1(x_1)=\cos(2\pi x_1), \quad V_2(x_2)=x_2^2/2, \quad V_3(x_3)=x_3^2/2.
\end{equation*}
The initial condition is given by
\begin{equation*}
\psi_0(\mathbf x) = 2^{-\frac{5}{2}} \pi^{-\frac{3}{4}} (x_1+\rmi x_2)\exp\left(-x_1^ 2/4 -x_2^ 2/4 -x_3^ 2/4 \right).
\end{equation*}

This equation could be integrated using any of the iterative methods considered in the previous section. However, for reasons of efficiency a time splitting approach is commonly employed. This treats the Laplacian and the potential part of the equations separately. For the former the fast Fourier transform (FFT) can be employed, while an analytic solution is available for the latter. The two partial flows are then combined by means of the Strang splitting scheme. For more details on this Time Splitting Fourier Pseudospectral method (TSFP) we refer the reader to \cite{JMS11}.

Another approach is to use a Hermite pseudospectral space discretization. This has the advantage that harmonic potentials are treated exactly, which is desirable in many applications. However, for most of the other potentials, the resulting matrices are full which, for traditional integration schemes, means that using a Hermite pseudospectral discretization is not competitive with respect to TSFP. However, as long as the potential is in Kronecker form, we can employ the $\mu$-mode integrator to perform computations very efficiently. Moreover, the resulting method based on the $\mu$-mode integrator combined with a Hermite pseudospectral space discretization can take arbitrarily large time steps without incurring any time discretization error (as it is exact in time). We call this scheme the Hermite Kronecker Pseudospectral method (HKP).

Before proceeding, let us note that for the TSFP method it is necessary to truncate the unbounded domain. In order to relate the size of the truncated domain to the chosen degrees of freedom, we considered that, in practice, in the HKP method the domain is implicitly truncated. This truncation is given by the convex hull of the quadrature points necessary to compute the Hermite coefficients corresponding to the initial solution. For any choice of degrees of freedom of the TSFP method, we decided to truncate the unbounded domain to the corresponding convex hull of the quadrature points of the HKP method. In this way, for the same degrees of freedom, the two methods use the same amount of information coming from the same computational domain.

\begin{figure}
    \centering
%
%
\definecolor{mycolor1}{rgb}{0.00000,0.44700,0.74100}%
\definecolor{mycolor2}{rgb}{0.85000,0.32500,0.09800}%
\begin{tikzpicture}

\begin{axis}[%
width=4.163in,
height=2.236in,
at={(0.772in,0.487in)},
scale only axis,
xmode=log,
xmin=1e-04,
xmax=5e4,
xminorticks=true,
xlabel style={font=\color{white!15!black}},
xlabel={Wall-clock time (s)},
xlabel style={font=\footnotesize},
xticklabel style={font=\footnotesize},
ymode=log,
ymin=4.4e-05,
ymax=1.3e-01,
yminorticks=true,
ylabel style={font=\color{white!15!black}},
ylabel style={font=\footnotesize},
ylabel={Accuracy},
yticklabel style={font=\footnotesize},
axis background/.style={fill=white},
legend style={legend cell align=left, align=left, draw=white!15!black},
legend style = {font=\footnotesize}
]

\addplot [thick, color=mycolor2, mark=asterisk, mark options={solid, mycolor2, scale = 1.5}]
table[row sep=crcr]{%
    0.00454	0.07\\
    0.047031	0.007\\
    0.165315	0.0007\\
    0.492848	7e-05\\
};
\addlegendentry{HK}

\addplot [thick, color=mycolor1, mark=o, mark options={solid, mycolor1, scale = 1.5}]
  table[row sep=crcr]{%
0.039802	0.07\\
1.276825	0.007\\
6.984421	0.0007\\
69.685792	7e-05\\
};
\addlegendentry{TSFP}

\node[left, align=left, font = \footnotesize]
at (axis cs:0.004,0.07) {$N=40^3$};
\node[left, align=left, font = \footnotesize]
at (axis cs:0.04,0.007) {$N=80^3$};
\node[left, align=left, font = \footnotesize]
at (axis cs:0.13,0.0007) {$N=120^3$};
\node[left, align=left, font = \footnotesize]
at (axis cs:0.4,0.00007) {$N=160^3$};

\node[right, align=left, font = \footnotesize]
at (axis cs:0.04,0.07) {$N=40^3$, $\text{ts}=8$};
\node[right, align=left, font = \footnotesize]
at (axis cs:1.5,0.007) {$N=80^3$, $\text{ts}=32$};
\node[right, align=left, font = \footnotesize]
at (axis cs:8.5,0.0007) {$N=120^3$, $\text{ts}=64$};
\node[right, align=left, font = \footnotesize]
at (axis cs:90.4,0.00007) {$N=160^3$, $\text{ts}=256$};

\end{axis}

\end{tikzpicture}%
    \caption{Precision diagram for the integration of the Schr\"odinger equation with a time independent potential (\ref{eq: schtind}) up to $T=1$. The number of degrees of freedom $N$ and the number of time steps (ts) are varied in order to achieve a result which is accurate up to the given tolerance. The reference solution has been computed by the HKP method with $N=300^3$. }
    \label{fig: SHTIND1}
\end{figure}
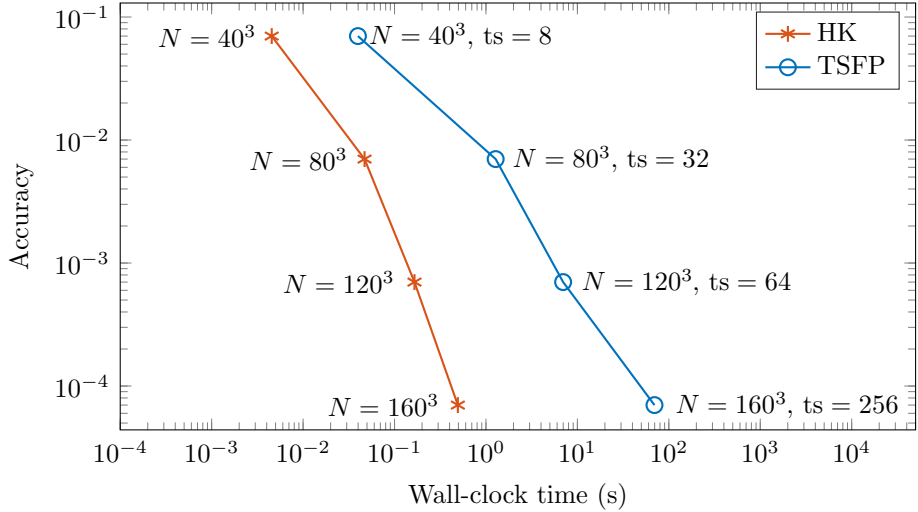

The TSFP and the HKP methods are compared in Figure \ref{fig: SHTIND1}. In both cases, we consider a constant number of space discretization points $n_\mu=n$ for every direction $\mu = 1,2,3$ (total number of degrees of freedom $N=n^3$) and integrate the equation until final time $T=1$ with constant time step size. We see that in terms of wall-clock time the HK method outperforms the TSFP scheme for all levels of accuracy considered here. Also note that the difference in performance increases as we move to more stringent tolerances. The reason for this is that the splitting error forces the TSFP scheme to take relatively small time steps.

\subsection{Schr\"odinger equation with time dependent potential \label{sec:schrtd}}

Let us now consider the Schr\"odinger equation
\begin{equation}
\left\{
\begin{aligned}
\partial_t \psi(t,\mathbf{x})&= H(t,\mathbf{x})\psi(t,\mathbf{x}), \qquad \mathbf{x}\in\RR^3 ,\quad t\in[0,1] \\[3pt]
\psi(0,\mathbf{x}) &= 2^{-\frac{5}{2}} \pi^{-\frac{3}{4}} (x_1+\rmi x_2)\exp\left(-x_1^ 2/4 -x_2^ 2/4 -x_3^ 2/4 \right),
\end{aligned}
\right.
\label{eq:sch}
\end{equation}
where the Hamiltonian is given by
\begin{equation*}
H(\mathbf{x},t) = \frac{\mathrm{i}}{2}\Big(\Delta -x_1^2-x_2^2-x_3^2 - 2x_3 \sin^2t\Big).
\end{equation*}
Note that the potential is now time dependent, as opposed to the case presented in section \ref{sec:schrti}. Such potentials commonly occur in applications, e.g.~when studying laser-atom interactions (see, for example, \cite{PKM94}).

Similarly to what we did in the time independent case, we can use a time splitting approach: the Laplacian part can still be computed efficiently in Fourier space, but now the potential part has no known analytical solution. Hence, for the numerical solution of the latter, we will employ an order two Magnus integrator, also known as the exponential midpoint rule. Let
\begin{equation*}
u'(t) = A(t)u(t)
\end{equation*}
be the considered ODE with time dependent coefficients, and let $u_n$ be the numerical approximation to the solution at time $t_n$. Then, the exponential midpoint rule provides the numerical solution
\begin{equation} \label{eq:mag2}
u_{n+1} = \exp\bigl(\kk_n A(t_n+\kk_n/2)\bigr)u_n
\end{equation}
at time $t_{n+1} = t_n + \kk_n$, where $\kk_n$ denotes the chosen step size.
The two partial flows are then combined together by means of the Strang splitting scheme. We call this scheme the Time Splitting Fourier Magnus Pseudospectral method (TSFMP). For the domain truncation needed in this approach, the same reasoning as in the time independent case applies.

Another technique is to perform a Hermite pseudospectral space discretization. However, as opposed to the case in section \ref{sec:schrti}, the resulting ODE  cannot be integrated exactly in time. For the time discretization, we will then use the order two Magnus integrator (\ref{eq:mag2}). We call the resulting scheme Hermite Kronecker Magnus Pseudospectral method (HKMP).

The results of the experiments are depicted in Figure \ref{fig: sch_td}. In both cases, we consider a constant number of space discretization points $n_\mu=n$ for every direction $\mu = 1,2,3$ (total number of degrees of freedom $N=n^3$) and solve the equation until final time $T=1$ with constant time step size. Moreover, concerning the TSFMP method, we integrate the subflow corresponding to the potential part with a single time step. Again, as we observed in the time independent case, the HKMP method outperforms the TSFMP scheme in any case. Notice in particular that, for the chosen degrees of freedom and time steps, the TSFMP method is not able to reach an accuracy of 1e-07, while the HKMP is.

\begin{figure}
    \centering
%
%
\definecolor{mycolor1}{rgb}{0.00000,0.44700,0.74100}%
\definecolor{mycolor2}{rgb}{0.85000,0.32500,0.09800}%
\begin{tikzpicture}

\begin{axis}[%
width=4.263in,
height=2.236in,
scale only axis,
xmode=log,
xmin=0.0000238496389381291,
xmax=500.8496389381292,
xminorticks=true,
xlabel style={font=\color{white!15!black}},
xlabel={Wall-clock time (s)},
xlabel style={font=\footnotesize},
xticklabel style={font=\footnotesize},
ymode=log,
ymin=3.26885973441867e-08,
ymax=0.0302668130745556,
yminorticks=true,
ylabel style={font=\color{white!15!black}},
ylabel={Accuracy},
ylabel style={font=\footnotesize},
yticklabel style={font=\footnotesize},
axis background/.style={fill=white},
legend style={legend cell align=left, align=left, draw=white!15!black},
legend style = {font=\footnotesize}
]

\addplot [thick, color=mycolor2, mark=asterisk, mark options={solid, mycolor2, scale = 1.5}]
table[row sep=crcr]{%
    0.002093	0.01\\
    0.006021	0.001\\
    0.029926	0.0001\\
    0.061763	1e-05\\
    0.248802	1e-06\\
    1.043454	1.02e-07\\
};
\addlegendentry{HKMP}

\addplot [thick, color=mycolor1, mark=o, mark options={solid, mycolor1, scale = 1.5}]
  table[row sep=crcr]{%
0.007921	0.01\\
0.082053	0.001\\
0.239208	0.0001\\
0.834867	1e-05\\
3.364443	1e-06\\
};
\addlegendentry{TSFMP}

\node[left, align=left, font = \footnotesize]
at (axis cs:0.002,9e-03) {$N=20^3$, $\text{ts}=2$};
\node[left, align=left, font = \footnotesize]
at (axis cs:0.005,9e-04) {$N=20^3$, $\text{ts}=8$};
\node[left, align=left, font = \footnotesize]
at (axis cs:0.025,9e-05) {$N=20^3$, $\text{ts}=32$};
\node[left, align=left, font = \footnotesize]
at (axis cs:0.055,9e-06) {$N= 30^3$, $\text{ts} = 64$};
\node[left, align=left, font = \footnotesize]
at (axis cs:0.2,9e-07) {$N = 30^3$, $\text{ts} = 256$};
\node[left, align=left,  font = \footnotesize]
at (axis cs:0.8,9e-08) {$N = 40^3$, $\text{ts} = 512$};

\node[right, align=left, font = \footnotesize]
at (axis cs:0.012,9e-03) {$N=20^3$, $\text{ts}=4$};
\node[right, align=left, font = \footnotesize]
at (axis cs:0.1,9e-04) {$N=30^3$, $\text{ts}=16$};
\node[right, align=left, font = \footnotesize]
at (axis cs:0.3,9e-05) {$N=30^3$, $\text{ts}=64$};
\node[right, align=left, font = \footnotesize]
at (axis cs:0.9,9e-06) {$N = 40^3$, $\text{ts} = 128$};
\node[right, align=left, font = \footnotesize]
at (axis cs:3.5,9e-07) {$N = 40^3$, $\text{ts} = 512$};

\end{axis}

\begin{axis}[%
width=4.211in,
height=1.105in,
at={(0in,0in)},
scale only axis,
xmin=0,
xmax=1,
ymin=0,
ymax=1,
axis line style={draw=none},
ticks=none,
axis x line*=bottom,
axis y line*=left
]
\end{axis}
\end{tikzpicture}%
    \caption{Precision diagram for the integration of the Schr\"odinger equation with a time dependent potential (\ref{eq:sch}) up to $T=1$. The number of degrees of freedom $N$ and the number of time steps (ts) are varied in order to achieve a result which is accurate up to the given tolerance. The reference solution has been computed by the HKMP method with $N=100^3$ and $\text{ts}=2048$.}
    \label{fig: sch_td}
\end{figure}

\subsection{Nonlinear Schr\"odinger/Gross--Pitaevskii equation \label{sec: gpe}}

In this section we consider the nonlinear Schr\"odinger equation
\begin{equation} \label{eq:vortex}
\partial_t \psi(t,\mathbf{x})=\frac{\rmi}{2}\Delta \psi(t,\mathbf{x})+ \frac{\rmi}{2}\Big(1-\lvert\psi(t,\mathbf{x})\rvert^2\Big)\psi(t,\mathbf{x}),
\end{equation}
which is also known as Gross--Pitaevskii equation.  The unknown  $\psi$ represents the wave function, $\mathbf x\in\RR^3$, $t\in[0,25]$, and the initial condition is constituted by the superimposition of two straight vortices in a background density $\lvert \psi_\infty\rvert^2=1$, in order to replicate the classical experiment of vortex reconnection (see \cite{CZ18} and the references therein for more details).

The initial datum and the boundary conditions given by the background density make it  quite difficult to use artificial periodic boundary conditions in a truncated domain, unless an expensive mirroring of the domain in the three dimensions is carried out.
Therefore, in order to solve (\ref{eq:vortex}) numerically, we consider the Time Splitting Finite Difference method proposed in \cite{CZ18}. More specifically, we truncate the unbounded domain to $\mathbf x \in[-20,20]^3$ and discretize by non-uniform finite differences with homogeneous Neumann boundary
conditions. The number $n_\mu$ of discretization points is the same in each direction, i.e. $n_\mu=n$, with $\mu=1,2,3$. After a proper transformation of variables in order to recover symmetry, we end up with a system of ODEs of the form
\begin{equation*}
\boldsymbol{\psi}'(t)=\frac{\rmi}{2}M_W\boldsymbol{\psi}(t)+
\frac{\rmi}{2}\Bigl(1-W^{-1}\lvert \boldsymbol{\psi}(t)\rvert^2\Bigr)\boldsymbol{\psi}(t),
\end{equation*}
where $M_W$ is a matrix in Kronecker form and $W$ is a diagonal weight matrix. Then, we employ a Strang splitting scheme for the time integration, in which the linear part is solved either by means of the $\mu$-mode integrator or by using the iterative methods indicated at the beginning of section \ref{sec:numerical}. The nonlinear subflow is integrated exactly.

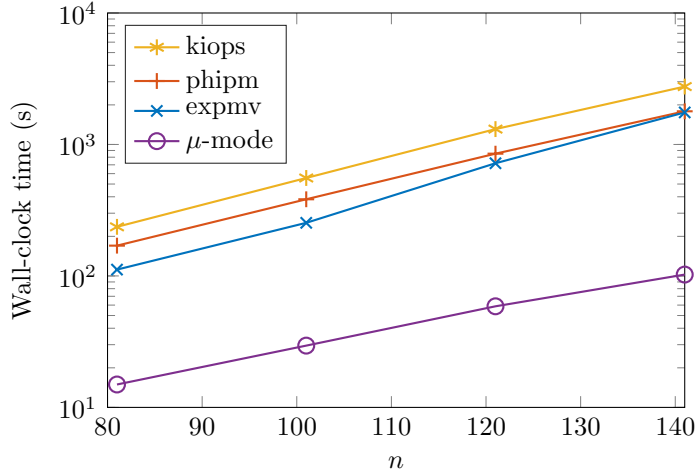
\begin{figure}
    \centering
%
%
\definecolor{mycolor1}{rgb}{0.92941,0.69412,0.12549}%
\definecolor{mycolor2}{rgb}{0.85098,0.32549,0.09804}%
\definecolor{mycolor3}{rgb}{0.00000,0.44706,0.74118}%
\definecolor{mycolor4}{rgb}{0.49412,0.18431,0.55686}%
\begin{tikzpicture}

\begin{axis}[%
width=3.021in,
height=2.066in,
at={(0.758in,0.481in)},
scale only axis,
xmin=80,
xmax=141,
xlabel style={font=\color{white!15!black}},
xlabel style={font=\footnotesize},
xticklabel style={font=\footnotesize},
xlabel={$n$},
ymode=log,
ymin=10,
ymax=10000,
yminorticks=true,
ylabel style={font=\color{white!15!black}},
ylabel={Wall-clock time (s)},
ylabel style={font=\footnotesize},
yticklabel style={font=\footnotesize},
axis background/.style={fill=white},
legend style={at={(0.03,0.97)}, anchor=north west, legend cell align=left, align=left, draw=white!15!black},
legend style={font=\footnotesize}
]
\addplot [thick, color=mycolor1, mark=asterisk, mark options={solid, mycolor1, scale = 1.5}]
  table[row sep=crcr]{%
81	236.271001\\
101	555.944092\\
121	1308.657725\\
141	2756.812657\\
};
\addlegendentry{kiops}

\addplot [thick, color=mycolor2, mark=+, mark options={solid, mycolor2, scale = 1.5}]
  table[row sep=crcr]{%
81	169.98812\\
101	382.918616\\
121	851.047256\\
141	1787.716756\\
};
\addlegendentry{phipm}

\addplot [thick, color=mycolor3, mark=x, mark options={solid, mycolor3, scale = 1.5}]
  table[row sep=crcr]{%
81	111.717933\\
101	253.489479\\
121	719.591097\\
141	1752.660457\\
};
\addlegendentry{expmv}

\addplot [thick, color=mycolor4, mark=o, mark options={solid, mycolor4, scale = 1.5}]
  table[row sep=crcr]{%
81	14.936964\\
101	29.509663\\
121	58.744092\\
141	102.402676\\
};
\addlegendentry{$\mu\text{-mode}$}

\end{axis}
\end{tikzpicture}%
    \caption{Wall-clock time (in seconds) for the integration of (\ref{eq:vortex}) up to $T=25$ as a function of $n$ (total number of degrees of freedom $N=n^3$). A constant time step size $\kk=0.1$ is employed.}
    \label{fig: vortex}
\end{figure}

The results of the experiment are presented in Figure \ref{fig: vortex}. The $\mu$-mode integrator outperforms \texttt{expmv} by approximately a factor of $7$. The speedup compared to both \texttt{phipm} and \texttt{kiops} is even larger.

\section{Implementation on multi-core CPUs and GPUs \label{sec:implementation}}

It has increasingly been realized that in order to fully exploit present and future high-performance computing systems  we require algorithms that parallelize well and which can be implemented efficiently on accelerators, such as GPUs \cite{ascac}. In particular, for GPU computing much research effort has been undertaken to obtain efficient implementations (see, e.g., \cite{auer2018magnus,PIConGPU2013,einkemmer2016mixed,einkemmer2017evaluation,einkemmer2020semi,lyakh2015,mehrenberger2013vlasov,rawat2018,sandroos2013multi,wiesenberger2019reproducibility}).

In this section we will consider an efficient implementation of the proposed $\mu$-mode integrator on multi-core CPUs and GPUs. We note that all modern hardware platforms are much better at performing floating point operations (such as addition and multiplication) than they are at accessing data in memory. This favors algorithms with a high flop/byte ratio; that is, algorithms that perform many floating point operations for every byte that is loaded from or written to memory. The $\mu$-mode product of a square matrix for an array of size $n_1 \times \dots \times n_{\mu-1} \times n_\mu \times n_{\mu+1} \times \dots \times n_d$ is computed using a matrix-matrix multiplication of size $n_\mu \times n_\mu$ times $n_\mu \times (n_1 \cdots n_{\mu-1} n_{\mu+1} \cdots n_{d})$, see section \ref{sec:mumod} for more details. For moderate $n_\mu$ the relatively small $n_\mu \times n_\mu$ matrix can be kept in cache and thus $\mathcal{O}(n_\mu N)$ arithmetic operations are performed compared to $\mathcal{O}(N)$ memory operations, where $N=n_1 \cdots n_d$ is the total number of degrees of freedom. Thus, the flop/byte ratio of the algorithm is $\mathcal{O}(n_\mu)$, which makes it ideally suited to modern computer hardware. This is particularly true when the $\mu$-mode integrator is compared to an implicit scheme implemented with sparse matrix-vector products. In this case the flop/byte ratio is only $\mathcal{O}(1)$ and modern CPU and GPUs will spend most of their time waiting for data that is fetched from memory.

To make this analysis more precise, we have to compare the flop/byte ratio of the algorithm to that of the hardware. For the benchmarks in this section we will use a multi-core CPU system based on a dual socket Intel Xeon Gold 5118 with $2 \times 12$ cores. The system has a peak floating point performance of $1.8$ TFlops/s (double precision) and a theoretical peak memory bandwidth of $256$ GB/s. Thus, during the time a double precision floating point number is fetched from memory approximately $56$ arithmetic operations can be performed. In addition, we will use a NVIDIA V100 GPU with 7.5 TFlop/s double precision performance and 900 GB/s peak memory bandwidth (approximately $67$ arithmetic operations can be performed for each number that is fetched from memory). Due to their large floating point performance we expect the algorithm to perform well on GPUs. A feature of the V100 GPU is that it contains so-called tensor cores that can dramatically accelerate half-precision computations (up to 125 Tflops/s). Tensor cores are primarily designed for machine learning tasks, but they can also be exploited for matrix-matrix products (see, e.g., \cite{abdelfattah2019,markidis2018}).

For reasonably large $n_\mu$ the proposed $\mu$-mode integrator is thus compute bound. However, since very efficient (close to the theoretical peak performance) matrix-matrix routines are available on both of these platforms, one can not be entirely indifferent towards memory operations. There are two basic ways to implement the algorithm. The first is to explicitly form the $n_\mu \times (n_1 \cdots n_{\mu-1} n_{\mu+1} \cdots n_{d})$ matrix.
This has the advantage that a single matrix-matrix multiplication (gemm) can be used to perform each $\mu$-mode product and that the corresponding operands have the proper sequential memory layout. The disadvantage is that a permute operation has to be performed before each $\mu$-mode product is computed. This is an extremely memory bound operation with strided access for which the floating point unit in the CPU or GPU lies entirely dormant. Thus, while this is clearly the favored approach in a \matlab{} implementation, it does not achieve optimal performance. The approach we have chosen in this section is to directly perform the $\mu$-mode products on the multi-dimensional array stored in memory (without altering the memory layout in between such operations).

Both Intel MKL and cuBLAS provide appropriate batched gemm routines (\texttt{cblas\_gemm\_batch} for Intel MKL and \texttt{cublasGemmStridedBatched} for cuBLAS) that are heavily optimized, and we will make use of those library functions in our implementation (for more details on these routines we refer to \cite{cecka2017}). Our code is written in C++ and uses CUDA for the GPU implementation. 

Before proceeding, let us briefly discuss how the $\mu$-mode integrator would perform in a distributed memory setting (i.e.~when parallelized using MPI). Since, in general, the matrix exponentials are full matrices, each degree of freedom along a coordinate axis couples with each other degree of freedom on that same axis. This data communication pattern is similar to computing a FFT. Thus, we would expect the $\mu$-mode product to scale comparable to FFT on a distributed memory system. This would be worse than a stencil code. However, one should keep in mind that the $\mu$-mode integrator can take much larger time steps. Thus, the overall communication overhead to compute the solution at a specified final time could still be larger for an explicit or an iterative method.

In the remainder of this section we will present benchmark results for our implementations.
The speedups are always calculated as ratio between the wall-clock time needed by the CPU and the one needed by the GPU.

\subsection{Heat equation}

We consider the same problem as in section \ref{sec:heat}, Test 1. The wall-clock time for computing the matrix exponentials and a single time step of the proposed algorithm is listed in Table \ref{tab:he}. We consider both a CPU implementation using MKL (double and single precision) and a GPU implementation based on cuBLAS (double, single, and half precision). The GPU implementation outperforms the CPU implementation by a factor of approximately $13$. Using half-precision computations on the GPU results in another performance increase by approximately a factor of $2$. The relative error with respect to the analytical solution reached by the double precision and single precision, for both CPU and GPU and the values of $n$ under consideration, are 8.22e-05, 3.66e-05, 2.06e-05, 1.47e-05. Results in half precision are not reported as the accuracy of the method is lower than the precision itself.

\begin{table}[!ht]
    \centering
    \bgroup\def\arraystretch{1.3}
    \begin{tabular}{c||c||ccc||ccc||c}
        $n$ & exp & \multicolumn{3}{c||}{double} & \multicolumn{3}{c||}{single} & half \\
        \hline
        & & CPU & GPU & speedup & CPU & GPU & speedup & GPU \\
        200 & 2.92 & 38.39 & 2.66 & 14.4x & 19.48 & 1.33 & 14.6x & 0.39 \\
        300 & 4.88 & 136.17 & 8.90 & 15.3x & 81.65 & 5.27 & 15.5x & 2.73 \\
        400 & 10.14 & 310.11 & 29.88 & 10.4x & 161.97 & 16.89 & 9.6x & 6.68 \\
        500 & 17.74 & 711.07 & 52.86 & 13.5x & 373.36 & 30.51 & 12.2x & 15.43
    \end{tabular}
    \egroup
    \caption{Wall-clock time for the heat equation (\ref{eq:he}) discretized using second-order centered finite differences with $n^3$ degrees of freedom. The time for computing the matrix exponentials (exp) and for one step of the $\mu$-mode integrator are listed (in ms). The speedup is the ratio between the single step performed in CPU and GPU, in double and single precision. The matrix exponential is always computed in double precision.}\label{tab:he}
\end{table}

For a number of simulations conducted we observed a drastic reduction in performance for single precision computations when using Intel MKL. To illustrate this we consider the heat equation
\begin{equation} \label{eq:hedenorm}
\left\{
\begin{aligned}
\partial_t u(t,\mathbf{x}) &= \Delta u(t,\mathbf{x}) ,\qquad \mathbf{x}\in \left[-\tfrac{11}4,\tfrac{11}4\right]^3,\quad t\in[0,1], \\[5pt]
u(0,\mathbf{x}) &= \left(x_1^4+x_2^4+x_3^4\right)\exp\bigl(-x_1^4-x_2^4-x_3^4\bigr) 
\end{aligned}
\right.
\end{equation}
with (artificial) Dirichlet boundary conditions, discretized in space as above.
From Table \ref{tab:nodenorm} we see that the performance of single precision computations with Intel MKL can be worse by a factor of $3.5$ compared to double precision, which obviously completely defeats the purpose of doing so. The reason for this performance degradation are so-called denormal numbers, i.e.~floating point numbers with leading zeros in the mantissa. Since there is no reliable way to disable denormal numbers on modern x86-64 systems, we avoid them by scaling the initial value in an appropriate way. Since this is a linear problem, the scaling can easily be undone after the computation. The results with the scaling workaround, listed in Table \ref{tab:nodenorm}, now show the expected behavior (that is, single precision computations are approximately twice as fast as double precision ones). We note that this is \textit{not} an issue with our $\mu$-mode integrator but rather an issue with Intel MKL. The cuBLAS implementation is free from this artifact and thus no normalization is necessary on the GPU.

\begin{table}[!ht]
    \centering
    \bgroup\def\arraystretch{1.3}
    \begin{tabular}{c||c||cc||cc||c||c}
        $n$ & exp & \multicolumn{2}{c||}{double} & \multicolumn{2}{c||}{single} & scaled single & half \\
        \hline
        & & CPU & GPU  & CPU & GPU & CPU & GPU \\
        200 & 2.92 & 38.80 & 2.64  & 92.19 & 1.34 & 19.98 &  0.38 \\
        300 & 6.01 & 157.41 & 8.87  & 385.84 & 5.22 &  71.24 & 2.71 \\
        400 & 13.40 & 314.96 & 29.85  & 1059.78 & 16.86 & 154.84 & 6.67 \\
        500 & 30.19 & 702.48 & 52.92  & 2567.56 & 30.42 & 367.34 & 13.44
    \end{tabular}
    \egroup

    \caption{Wall-clock time for the heat equation (\ref{eq:hedenorm}) discretized using second order centered finite differences  with $n^3$ degrees of freedom. The performance degradation in CPU due to denormal numbers disappears when using the scaling workaround (scaled single). Speedups are not computed in this case.} \label{tab:nodenorm}
\end{table}

\subsection{Schr\"odinger equation with time independent potential}

We consider the Schr\"odinger equation with time independent potential from section \ref{sec:schrti}. The equation is integrated up to $T=1$ in a single step, as for this problem no error is introduced by the $\mu$-mode integrator. For the space discretization the Hermite pseudospectral discretization is used. The results for both the CPU and GPU implementation are listed in Table \ref{tab:schhtind}. The GPU implementation, for both single and double precision, shows a speedup of approximately $15$ compared to the CPU implementation.

\begin{table}[!ht]
    \centering
    \bgroup\def\arraystretch{1.5}
        \begin{tabular}{c||c|ccc||c|ccc}
            $n$ & \multicolumn{4}{c||}{double} & \multicolumn{4}{c}{single} \\
            \hline
            & exp & CPU  & GPU & speedup & exp & CPU  & GPU & speedup \\
            127 & 5.56 & 20.89 & 1.27 & 16.4x & 4.71 & 13.71 & 0.64 & 21.4x \\
            255 & 8.31 & 224.13 & 16.02 & 13.9x & 5.16 & 134.21 & 8.11 & 16.5x \\
            511 & 50.79 & 3121.42 & 219.13 & 14.2x & 28.01 & 1824.93 & 119.46 & 15.2x
        \end{tabular}
    \egroup
    \caption{Wall-clock time for the linear Schr\"odinger equation with time independent potential (\ref{eq: schtind}) integrated with the HKP method ($n^3$ degrees of freedom). The time for computing the matrix exponential (exp) and for one step of the $\mu$-mode integrator are listed (in ms). The speedup is the ratio between the single step performed in CPU and GPU, in double and single precision.}\label{tab:schhtind}
\end{table}

\subsection{Schr\"odinger equation with time dependent potential}

We consider once again the Schr\"odinger equation with the time dependent potential from section \ref{sec:schrtd} solved with the HKMP method. The equation is integrated up to $T=1$ with time step $\kk = 0.02$. The results are given in Table \ref{tab:schhtd}. In this case, the matrix exponential changes  as we evolve the system in time. Thus, the performance of computing the matrix exponential has to be considered alongside the $\mu$-mode products. On the CPU this is not an issue as the time required for the matrix exponential is significantly smaller than the time required for the $\mu$-mode products. However, for the GPU implementation and small problem sizes it is necessary to perform the matrix exponential on the GPU as well. To do this we have implemented an algorithm based on a Taylor backward stable approach. Overall, we observe a speedup of approximately $15$ by going from the CPU to the GPU (for both single and double precision).

\begin{table}[h!]
    \centering

     \bgroup\def\arraystretch{1.3}
    \begin{tabular}{c||c|c|c|c|c|c}
        $n$ & \multicolumn{6}{c}{double} \\
        \hline
        & exp (ext) & \multicolumn{2}{c|}{CPU} & \multicolumn{2}{c|}{GPU} & speedup \\
        & & exp (int) & $\mu$-mode & exp (int) & $\mu$-mode & \\
        127 & 0.02 & 2.56 & 19.38 & 0.37 & 1.05 & 15.3x \\
        255 & 0.05 & 4.52 & 200.46 & 0.66 & 13.79 & 14.2x \\
        511 & 0.07 & 29.71 & 3043.88 & 2.38 & 213.21 & 14.3x
    \end{tabular}
    \egroup

    \vspace{0.5cm}

    \bgroup\def\arraystretch{1.3}
        \begin{tabular}{c||c|c|c|c|c|c}
            $n$ & \multicolumn{6}{c}{single} \\
            \hline
            & exp (ext) & \multicolumn{2}{c|}{CPU} & \multicolumn{2}{c|}{GPU} & speedup \\
            & & exp (int) & $\mu$-mode & exp (int) & $\mu$-mode & \\
            127 & 0.01 & 2.16 & 12.51 & 0.25 & 0.54 & 18.9x \\
            255 & 0.03 & 2.88 & 100.35 & 0.34 & 7.01 & 13.9x \\
            511 & 0.05 & 14.25 & 1600.86 & 1.09 & 108.31 & 14.8x
        \end{tabular}
    \egroup
 
    \caption{Wall-clock time for the Schr\"odinger equation with time dependent potential (\ref{eq:sch}) integrated with the HKMP method ($n^3$ degrees of freedom). The time for computing the matrix exponentials and for one step of the $\mu$-mode integrator is  listed (in ms). The acronym exp (ext) refers to exponentiation of the time independent matrices, which are diagonal, while exp (int) refers to the time dependent ones that have to be computed at each time step. The speedup is the ratio between the single step performed in CPU and GPU, in double precision (top) and single precision (bottom).}\label{tab:schhtd}

\end{table}

\section{Conclusions}\label{sec:conclusions}

We have shown that with the proposed $\mu$-mode integrator we can make use of modern computer hardware to efficiently solve a number of partial differential equations. In particular, we have demonstrated that for Schr\"odinger equations the approach can outperform well established integrators in the literature by a significant margin. This was also possible thanks to the usage of the $\mu$-mode product to efficiently compute spectral transforms, which can be beneficial even in applications that are not related to solving partial differential equations. The proposed integrator is particularly efficient on GPUs too, as we have demonstrated, which is a significant asset for running simulation on the current and next generation of supercomputers.

\section*{Acknowledgments}
The authors would like to thank the Italian Ministry of Instruction,
University and Research (MIUR) for partially supporting this research
with funds coming from PRIN Project 2017 (No. 2017KKJP4X entitled
``Innovative numerical methods for evolutionary partial differential equations
and applications'').

The authors acknowledge partial support from the Program Ricerca di Base
2019 of the University of Verona
entitled
``Geometric Evolution of Multi Agent Systems''.

FZ has received funding from the 
European Research Council (ERC) 
under the European Union’s Horizon 2020 
research and innovation programme (grant agreement No. 850941).

This work is supported by the Austrian Science Fund (FWF) -- project id: P32143-N32.

\bibliographystyle{plain}
\bibliography{kronbiblio}

\end{document}